\documentclass{amsart}
\usepackage{amsmath,amsfonts}
\usepackage{latexsym}
\title[K-theory of semi-periodic pseudos]
{K-theory of pseudodifferential operators
with semi-periodic symbols}
\author[S. T. Melo \and C. C. Silva]{Severino T. Melo 
\and C\'intia C. Silva}
\date{}
\newtheorem{thm}{Theorem}
\newtheorem{pro}{Proposition}
\newtheorem{lem}{Lemma}
\newtheorem{cor}{Corollary}

\begin{document}
\newcommand{\U}{\mbox{{\tt u}}}
\newcommand{\V}{\mbox{{\tt v}}}
\newcommand{\s}{{\mathbb S}}
\newcommand{\R}{{\mathbb R}}
\newcommand{\C}{{\mathbb C}}
\newcommand{\N}{{\mathbb N}}
\newcommand{\Z}{{\mathbb Z}}
\newcommand{\op}{operator}
\newcommand{\ops}{operators}
\newcommand{\psd}{pseudo\-dif\-fer\-en\-tial}
\newcommand{\Psd}{Pseudo\-dif\-fer\-en\-tial}
\newcommand{\cqd}{\hfill$\Box$}
\newcommand{\dx}{d\!x}
\newcommand{\dy}{d\!y}
\newcommand{\lt}{L^2({\mathbb R})}
\newcommand{\as}{{\mathcal A}}
\newcommand{\ac}{{\mathfrak A}}
\newcommand{\jc}{{\mathfrak J}}
\newcommand{\kc}{{\mathfrak K}}
\newcommand{\ec}{{\mathfrak E}}
\newcommand{\bc}{{\mathfrak B}}
\newcommand{\Cc}{{\mathfrak C}}
\newcommand{\img}{\text{\tt Im}\,\gamma}
\newcommand{\pf}{{\em Proof}: }
\newcommand{\cst}{C$^*$}
\newcommand{\ind}{\mbox{{\tt ind}}\,}

\begin{abstract} Let $\ac$ denote the \cst-algebra of 
bounded \ops\ on $\lt$ generated by: (i) all multiplications 
$a(M)$ by functions $a\in C[-\infty,+\infty]$, (ii) all multiplications by 
$2\pi$-periodic continuous functions, and (iii) all \ops\ of the form 
$F^{-1}b(M)F$, where $F$ denotes the Fourier transform and 
$b\in C[-\infty,+\infty]$. 
A given $A\in\ac$ is a Fredholm \op\ if and only if $\sigma(A)$ and 
$\gamma(A)$ are invertible, where $\sigma$
denotes the continuous extension of the usual principal symbol, while 
$\gamma$ denotes an operator-valued ``boundary principal symbol'' (the
``boundary'' here consists of two copies of the circle, one at each end of the
real line). We give two proofs of the fact that $K_0(\ac)$ is isomorphic to 
$\Z$ and that $K_1(\ac)$ is isomorphic to $\Z\oplus\Z$. We do it first by 
computing the connecting mappings in the six-term exact sequence associated to 
$\sigma$. For the second proof, we show that the image of $\gamma$ is 
isomorphic to the direct sum of two copies of the crossed product 
$C[-\infty,+\infty] \times_{\alpha}\Z$, where $\alpha$ denotes the 
translation-by-one automorphism. Its K-theory can be computed using the 
Pimsner-Voiculescu exact sequence, and that information suffices for the 
analysis of the standard cyclic exact sequence associated to $\gamma$. 

\end{abstract}
\maketitle

\begin{center}
\end{center}

\section*{Introduction}

Let $\ac$ denote the smallest \cst-subalgebra of the algebra of all bounded 
\ops\ on $\lt$ containing: 

\begin{enumerate}
\label{gen}

\item every multiplication $a(M)$ by an $a\in C[-\infty,+\infty]$, where
$[-\infty,+\infty]$ denotes the two-point compactification of $\R$ (i.e., 
$a$ is continuous on $\R$ and has limits at $+\infty$ and at $-\infty$), 

\item every multiplication by a $2\pi$-periodic continuous function, and 

\item every {\em Fourier multiplier} of the form $b(D)=F^{-1}b(M)F$, where $F$ 
denotes the Fourier transform and $b\in C[-\infty,+\infty]$.

\end{enumerate}
The structure of $\ac$ is described in the following four theorems,
which are rephrased versions of results first appearing in \cite{CM,M1},
where this algebra is regarded as a particular {\em comparison algebra};
i.e., as a member of a certain class, defined by Cordes \cite{C2}, of \cst-algebras generated
by pseudodifferential operators. These results are partially reproven in the first chapter 
of \cite{S}, avoiding reference to general results on comparison algebras. 

\begin{thm}\label{th1}Let $X$ denote the subset of 
$[-\infty,+\infty]\times \s^1$ of all points $(x,e^{i\theta})$ such that 
$|x|=\infty$ or $x=\theta\in\R$.
There exists a surjective \cst-algebra homomorphism 
\[
\begin{array}{rcl}
\sigma:\ac&\longrightarrow&C(X\times\{-\infty,+\infty\})\\
A&\longmapsto&\sigma_{_{A}}
\end{array}
\] 
such that, for each $m=(x,e^{i\theta})\in X$,
we have: 

$(\!$i$)$ $\sigma_{_{a(M)}}(m,+\infty)=\sigma_{_{a(M)}}(m,-\infty)=a(x)$ if 
$a\in C[-\infty,+\infty]$, 

$(\!$ii$)$ $\sigma_{_{a(M)}}(m,+\infty)=\sigma_{_{a(M)}}(m,-\infty)=a(\theta)$ 
if $a$ is a continuous $2\pi$-periodic function on $\R$, and 

$(\!$iii$)$ $\sigma_{_{b(D)}}(m,+\infty)=b(+\infty)$ and
$\sigma_{_{b(D)}}(m,-\infty)=b(-\infty)$ if 
$b\in C[-\infty,+\infty]$. 
\end{thm}

\begin{cor}\label{cor1}
Let $\ec$ denote the kernel of $\sigma$. The mapping
\begin{equation}\label{a/e}
\begin{array}{rcl}
\ac/\ec&\longrightarrow&C(X\times\{-\infty,+\infty\})\\
{[A]_\ec}&\longmapsto &\sigma_{_{A}}
\end{array}
\end{equation}
is a \cst-algebra isomorphism.
\end{cor}

For each $\varphi\in\R$, let $U_\varphi$ denote the operator on 
$L^2(\s^1)$ given by $U_\varphi f(z)=z^{-\varphi}f(z)$, $z\in\s^1$, and let
$Y_\varphi$ denote the operator on $\ell^2(\Z)$ obtained by
conjugating $U_\varphi$ with the discrete Fourier transform. Then,
$Y_\varphi$ is a smooth family of unitary operators such that, for all 
$j\in\Z$, $Y_j((u_k)_{k\in\Z})=(u_{j+k})_{k\in\Z}$.

Below we denote by $S^*\s^1=\s^1\times\{-1,+1\}$ the co-sphere bundle over the unit
circle $\s^1$, by $\bc$ the algebra of all bounded operators on 
$\ell^2(\Z)$, and the identity operator on $\ell^2(\Z)$ or on $L^2(\R)$ by $I$.

\begin{thm}\label{th2} 
There exists a \cst-algebra homomorphism 
\[
\begin{array}{rcl}
\gamma:\ac&\longrightarrow&C(S^*\s^1,\bc)\\
A&\longmapsto&\gamma_{_{A}}
\end{array}
\]
such that, for each $v=(e^{2\pi i\varphi},\pm 1)\in S^*\s^1$, we
have: 

$(\!$i$)$ $\gamma_{_{a(M)}}(v)=a(\pm\infty)I$ if 
$a\in C[-\infty,+\infty]$,

$(\!$ii$)$ $\gamma_{_{e_{k}(M)}}(v)=Y_{-k}$, for $e_k(x)=e^{ikx}$, 
$k\in\Z$, and 

$(\!$iii$)$ $\gamma_{_{b(D)}}(v)=Y_\varphi M^b_\varphi Y_{-\varphi}$, if 
$b\in C[-\infty,+\infty]$, 
with $M^b_\varphi$ defined by  
$M^b_\varphi((a_j)_{j\in\Z})=(b(j-\varphi)a_j)_{j\in\Z}$. 
\end{thm}

It is easy to show, and this is part of the statement of Theorem~\ref{th2}, that the function
$\R\ni\varphi\mapsto Y_\varphi M^b_\varphi Y_{-\varphi}\in\bc$ is 1-periodic if $b\in C[-\infty,+\infty]$.

Let $\as^\sharp$ denote the algebra of bounded functions on $\R$ generated
by $C[-\infty,+\infty]$ and by all $2\pi$-periodic continuous functions. 
It is closed. Let $C_0(\R)$ denote the algebra of all continuous functions
with zero limits at $+\infty$ and $-\infty$.

Let $\kc_\R$ denote the ideal of compact operators on $L^2(\R)$ and
$\kc_\Z$ denote the ideal of compact operators on $\ell^2(\Z)$. 

\begin{thm}
\label{th3} 
A dense *-subalgebra of the kernel of $\sigma$ 
 is given by the linear span of all operators 
of the form $a(M)b(D)+K$, with $a\in\as^\sharp$, $b\in C_0(\R)$, and 
$K\in\kc_\R$. The image of the restriction of $\gamma$ to $\ec$ is equal to 
$C(S^*\s^1,\kc_\Z)$ and the mapping 
\begin{equation}
\label{e/k}
\begin{array}{rcl}
\ec/\kc_\R&\longrightarrow&C(S^*\s^1,\kc_\Z)\\
{[A]_ {\kc_{\R}}}&\longmapsto&\gamma_{_{A}}
\end{array}
\end{equation}
is a \cst-algebra isomorphism. 
\end{thm}

\begin{thm}
The intersection of the kernels of $\sigma$ and $\gamma$ is equal to $\kc_\R$. Moreover, a given $A\in\ac$ is a Fredholm operator 
if and only if $\sigma(A)$ and $\gamma(A)$ are invertible.
\end{thm}

A similar Fredholm criterion had been given earlier by Rabinovich \cite{R} for 
an algebra larger than this. The above stated results were later 
extended  \cite{M2} to a comparison algebra on a cylinder of the form $\R\times B$, 
where $B$ is a compact Riemannian manifold.

In this paper, we give two proofs that $K_0(\ac)\cong\Z$ and $K_1(\ac)\cong\Z\oplus\Z$ (Theorem~\ref{main}): 
by computing the connecting mappings 
in the standard K-theory cyclic exact sequences associated to the two \cst-algebra short exact
sequences induced by the principal symbol $\sigma$ and by the boundary principal symbol $\gamma$. 
The crucial step in the second computation is finding a useful characterization of the image of 
$\gamma$, $\img$, and computing its K-groups. We do that by showing that $\img$ is
isomorphic to the direct sum of two copies of the crossed product 
$C[-\infty,+\infty]\times_{\alpha}\Z$, where $\alpha$ denotes the 
translation-by-one automorphism. Its K-theory can then be computed using the 
Pimsner-Voiculescu exact sequence.

Many other \cst-algebras generated by zero-order \psd\ \ops\ 
(see \cite{C1,CP,Lt,L,M2,MNS,Mo}, for example, but this is a very incomplete list)
have basically the same structure as our $\ac$; i.e., a two-step composition series of the form 
\begin{equation}\label{cs}
\kc\subset\ec\subset\ac,
\end{equation}
with $\ac/\ec$ commutative and 
\begin{equation}\label{e/ka}
\ec/\kc \cong C_0(Y)\otimes\kc,
\end{equation}
for some locally compact Hausdorff space $Y$.
In \eqref{cs}, $\kc$ is the ideal of compact operators on the Hilbert
space where $\ac$ acts, while in \eqref{e/ka} $\kc$ 
should rather be regarded as 
the compact ideal on a Hilbert space $H_b$ somehow associated to a ``boundary'', 
which can be the actual boundary of a manifold with boundary in some cases.
In most examples, the Fredholm property for an arbitrary operator $A\in\ac$ is 
equivalent to the invertibility of $\sigma(A)$ and $\gamma(A)$; where $\sigma$ 
(the principal symbol) is the composition of the Gelfand mapping for $\ac/\ec$ 
with the canonical projection on the quotient, and $\gamma$ (the boundary 
principal symbol) is a \cst-algebra homomorphism 
\[
\gamma:\ac\longrightarrow C_0(Y)\otimes\bc
\]
which extends the composition of the isomorphism in \eqref{e/ka} with the 
quotient projection. Here $\bc$ denotes the algebra of bounded operators
on $H_b$.

Given an algebra $\ac$ with such a structure, one might ask what pieces of 
information are needed to compute its K-theory, how to explicitly describe generators 
for its K-groups. A detailed understanding of the example treated 
in this paper will hopefully offer some insight for this more abstract
question. 

A strategy parallel to that of Section \ref{img} was used by Melo, Nest and Schrohe in \cite{MNS}. In that paper,
finding convenient descriptions of the image and the kernel of the boundary principal symbol was the key for a successful 
analysis of the cyclic exact sequence associated to $\gamma$, which lead to the computation of 
the K-theory of Boutet de Monvel's algebra on a compact manifold whose boundary is non-empty and has
torsion-free K-groups. With the help of K-theory tools a little more sofisticated than the standard cyclic sequence,
but still relying on the same description of the image and kernel of $\gamma$, that torsion-free assumption was later removed 
by Melo, Schick and Schrohe \cite{MSS}. For the computation of the K-theory of the comparison algebra with periodic 
multiplications on a cylinder studied in \cite{M2}, also the approach of finding a good description of the image of the 
boundary principal symbol seems to be more promissing than trying to compute directly the connecting mappings in the 
cyclic exact sequence associated to the principal symbol. 

Using the language of groupoids, Monthubert and Nistor \cite{MoN} extended Atiyah and Singer's \cite{AS} definition of topological 
index for manifolds with corners and used it to compute the connecting mappings in the K-theory cyclic exact sequence associated to 
the principal symbol of Melrose's b-calculus. For the algebra of suspended pseudo-differential operators, also introduced by 
Melrose \cite{Me}, Moroianu \cite{Mo} found an interesting connection between the Fredholm index of Dirac operators and the exponential 
mapping (the index mapping vanishes in this case) in the cyclic sequence associated to the principal symbol.

We use twice in this paper the following index theorem due to Cordes, Herman and Power \cite{CH,P}.

\begin{thm}
\label{thm4} 
Let $\Cc$ denote the \cst-algebra of bounded operators on $L^2(\R)$ generated by all
operators of type $(1)$ or $(3)$ $($listed at the beginning of this paper$)$. Let ${\mathbb M}$ denote the 
boundary points of the square compactification of $\R^2$,
\[
{\mathbb M}\,=\,\{(x,\xi)\in[-\infty,+\infty]\times[-\infty,+\infty];\,|x|+|\xi|=\infty\}.
\]
The assignment $a(M)b(D)\mapsto a(x)b(\xi)$ extends to a surjective \cst-algebra homomorphism 
$\bar\sigma:\Cc\to C({\mathbb M})$ with kernel $\kc_\R$. A given $A\in\Cc$ is a Fredholm operator
if and only if $\bar\sigma(A)$ does not vanish at any point of ${\mathbb M}$; in which case, its Fredholm index is 
equal to the winding number of $\bar\sigma(A)$ $($if ${\mathbb M}$ is identified with the circle in 
the canonical way$)$. 
\end{thm}

%
%
%
%
%
%

\section{The principal-symbol exact sequence}
\label{principal}

Before analysing the K-theory cyclic exact sequence associated to 
\begin{equation}
\label{seq}
0\to \ec/\kc_\R\to\ac/\kc_\R\to\ac/\ec\to 0,
\end{equation} 
we need first to describe the K-groups of $C(X)$, since it follows immediately 
from Corollary~\ref{cor1} that $\ac/\ec\cong C(X)\oplus C(X)$.

Choose \label{choice}a smooth non-increasing function
$c:\R\to[0,1]$ such that $c(x)=1$ if $x\leq -1/5$ and $c(x)=0$ if $x\geq 1/5$,
and let $b=1-c$. Later it will be clear why we require $c$ and $b$ to be 
constant outside this narrow interval. At first it would suffice to suppose only 
that $c(-\infty)=b(+\infty)=1$. Define also 
\[
l(x)=\left\{\begin{array}{lr}
e^{ix}\ &\mbox{if}\ x\geq 0\\
1 \ &\mbox{if}\ x< 0
\end{array}\right.\ \mbox{and}\ \ 
\tilde{l}(x)=\left\{\begin{array}{lr}
1\ &\mbox{if}\ x\geq 0\\
e^{ix} \ &\mbox{if}\ x< 0\end{array}\right..
\]
We may regard $\R$ as an open dense subset of $X$ (the precise statement is: the mapping 
$x\mapsto (x,e^{ix})$ is a homeomorphism of $\R$ onto an open dense subset of $X$),
and then denote also by $l$ and $\tilde{l}$ their unique continuous extensions to $X$.

Given a projection $p$ in a \cst-algebra $A$, and a unitary $u$ in $A$ or in its unitization 
$A^+$, we denote by $[p]_0$ and $[u]_1$ the elements they represent, respectively, in $K_0(A)$ 
and in $K_1(A)$. 

Given $f$ and $g$ in $C(\s^1)$, we denote by $(f,g)$ the function on $\{-\infty,+\infty\}\times\s^1$ 
which is $f$ over $-\infty$ and $g$ over $+\infty$, and denote by $\mbox{{\tt z}}$ the
identity function on $\s^1$. Then, for example, $(1,\mbox{{\tt z}})$ denotes the function which is
constant and equal to $1$ on $\{-\infty\}\times\s^1$ and maps $(+\infty,z)$ to $z$ for all $z\in\s^1$.
Analogous notation will be used for operator-valued funcions on $S^*\s^1$ and also for functions on 
$X\times\{-\infty,+\infty\}$. 

\begin{pro} \label{prop1}
$K_0(C(X))=\Z[1]_0$ and $K_1(C(X))=\Z[l]_1\oplus\Z[\tilde l]_1$.\label{kx}
\end{pro}

\pf Let us consider the exact sequence 
\begin{equation}
\label{comm}
0\to C_0(\R)\to C(X)\to C(\{-\infty,+\infty\}\times\s^1)\to 0
\end{equation}
defined by the restriction mapping $C(X)\to C(\{-\infty,+\infty\}\times\s^1)$.
We know that: $K_0(C_0(\R))=0$, $K_1(C_0(\R))=\Z[e^{2\pi ib}]_1$, 
\[ K_0(C(\{-\infty,+\infty\}\times\s^1))=\Z[(1,0)]_0\oplus\Z[(0,1)]_0,\] and
\[ K_1(C(\{-\infty,+\infty\}\times\s^1))=\Z[(\mbox{{\tt z}},1)]_1\oplus\Z[(1,\mbox{{\tt z}})]_1.\]

The standard $K$-theory cyclic sequence associated to \eqref{comm} then becomes
\begin{equation}
\def\mapup#1{\Big\uparrow\rlap{$\vcenter{\hbox{$\scriptstyle#1$}}$}}
\def\mapdn#1{\Big\downarrow\rlap{$\vcenter{\hbox{$\scriptstyle#1$}}$}}
\label{cyc}
\begin{array}{ccccc}
0&\longrightarrow&K_0(C(X))&\longrightarrow &\Z[(1,0)]_0\oplus\Z[(0,1)]_0
\\&&&&
\\\mapup{\delta_1}&&&&\mapdn{\delta_0}
\\&&&&
\\\Z[(\mbox{{\tt z}},1)]_1\oplus\Z[(1,\mbox{{\tt z}})]_1&\longleftarrow&K_1(C(X))
&\longleftarrow& \Z[e^{2\pi ib}]_1.
\end{array}
\end{equation}
The upper-right horizontal arrow in (\ref{cyc}) maps 
$[1]_0$ to $[(1,1)]_0$, while the lower-left one maps $[l]_1$ to 
$[(1,\mbox{{\tt z}})]_1$ and $[\tilde l]_1$ to $[(\mbox{{\tt z}},1)]_1$. 
The restriction of $b$ to $\{-\infty,+\infty\}\times\s^1$ is $(0,1)$. It then follows
from \cite[12.2.2]{Ro} that $\delta_0([(0,1)]_0)=-[e^{2\pi ib}]_1$. Analogously we get 
$\delta_0([(1,0)]_0)=-[e^{2\pi ic}]_1$. Noticing that $[e^{2\pi ib}]_1=-[e^{2\pi ic}]_1$, the 
proof can now be finished by elementary algebraic arguments. \cqd

Defining $A_1=l(M)b(D)+c(D)$, $A_2=\tilde l(M)b(D)+c(D)$, $A_3=b(D)+l(M)c(D)$
and $A_4=b(D)+\tilde l(M)c(D)$, we have $\sigma_{_{A_{1}}}=(1,l)$,
$\sigma_{_{A_{2}}}=(1,\tilde l)$, $\sigma_{_{A_{3}}}=(l,1)$ and $\sigma_{_{A_{4}}}=(\tilde l,1)$. 
Denoting by $[B]_{\ec}$ the class in $\ac/\ec$ of a $B\in\ac$, we obtain from the isomorphism (\ref{a/e})
and from Proposition~\ref{kx}:
\begin{equation}
\label{k1ae}
K_1(\ac/\ec)=\Z[[A_1]_\ec]_1\oplus\Z[[A_2]_\ec]_1\oplus\Z[[A_3]_\ec]_1\oplus\Z[[A_4]_\ec]_1.
\end{equation}
Since $\sigma_{_{b(D)}}=(0,1)$ and $\sigma_{_{c(D)}}=(1,0)$, we also get:
\begin{equation}
\label{k0ae}
K_0(\ac/\ec)=\Z[[b(D)]_\ec]_0\oplus\Z[[c(D)]_\ec]_0.
\end{equation}

Next we focus on the K-groups of $\ec/\kc_\R$. The isomorphism (\ref{e/k}) and standard results in 
K-theory of C$^*$-algebras imply at once that $K_0(\ec/\kc_\R)$ and $K_1(\ec/\kc_\R)$ are both 
isomorphic to $\Z\oplus\Z$. More detailed information can be given if we consider the
exact sequence 
\begin{equation}
\label{cindida}
0\to S\kc_\Z\oplus S\kc_\Z\to C(S^*\s^1,\kc_\Z)\to\kc_\Z\oplus\kc_\Z\to 0
\end{equation}
($S$ denotes suspension) induced by the \cst-algebra homomorphism
\begin{equation}
\label{nu}
\begin{array}{rcl}
\nu: C(S^*\s^1,\kc_\Z)&\longrightarrow&\kc_\Z\oplus\kc_\Z\\
 (f,g)&\longmapsto &(f(1),g(1)).
\end{array}
\end{equation}

The fact that (\ref{cindida}) splits implies that $\nu$ induces a $K_0$-isomorphism, while the inclusion 
of $S\kc_\Z\oplus S\kc_\Z=\ker\nu$ into $C(S^*\s^1,\kc_\Z)$ induces a $K_1$-isomorphism.
We know that if $E$ is any (we will choose one later) rank-one projection on $\ell^2(\Z)$, then 
\begin{equation}\label{ktc}
K_0(\kc_\Z)=\Z[E]_0
\end{equation} 
and, by Bott periodicity \cite[11.1.2]{Ro}, $K_1(C(\s^1,\kc_\Z))=\Z[I+(\mbox{{\tt z}}-1)E]_1$. 
Let  $\U$ and $\V$ denote the unitaries in $(C(S^*\s^1,\kc_\Z))^+$ defined by
\[
\U=(I+(\mbox{{\tt z}}-1)E,I)\ \ \mbox{and}\ \ \V=(I,I+(\mbox{{\tt z}}-1)E).
\]
From $C(S^*\s^1,\kc_\Z)\cong C(\s^1,\kc_\Z)\oplus C(\s^1,\kc_\Z)$, it follows that:
\begin{equation}
\label{k1ek}K_1(\ec/\kc_\R)\cong K_1(C(S^*\s^1,\kc_\Z))=\Z[\U]_1\oplus\Z[\V]_1,
\end{equation}
with the isomorphism in (\ref{k1ek}) being induced by (\ref{e/k}).
It also follows from the above remarks that
\begin{equation}
\label{k0ek}
K_0(\ec/\kc_\R)\cong K_0(C(S^*\s^1,\kc_\Z))\cong K_0(\kc_\Z)\oplus K_0(\kc_\Z),
\end{equation}
where the second isomorphism is induced by $\nu$, and the first by (\ref{e/k}).

Substituting the isomorphism (\ref{e/k}) into (\ref{seq}), we get:
\begin{equation}
\label{seqc}
0\to C(S^*\s^1,\kc_\Z)\to\ac/\kc_\R\to\ac/\ec\to 0.
\end{equation} 
Substituting (\ref{k1ae}), (\ref{k0ae}), (\ref{k1ek}) and (\ref{k0ek}) into 
the cyclic exact sequence in K-theory associated to (\ref{seqc}) we get:
\begin{equation}
\label{sceq}
\def\mapup#1{\Big\uparrow\rlap{$\vcenter{\hbox{$\scriptstyle#1$}}$}}
\def\mapdn#1{\Big\downarrow\rlap{$\vcenter{\hbox{$\scriptstyle#1$}}$}}
\begin{array}{ccccc}
K_0(\kc_\Z)\oplus K_0(\kc_\Z)
&\longrightarrow&K_0(\ac/\kc_\R)&\longrightarrow &\Z[[b(D)]_\ec]_0\oplus\Z[[c(D)]_\ec]_0
\\&&&&
\\\mapup{\delta_1}&&&&\mapdn{\delta_0}
\\&&&&
\\\Z[[A_1]_\ec]_1\oplus\cdots\oplus\Z[[A_4]_\ec]_1
&\longleftarrow&K_1(\ac/\kc_\R)&\longleftarrow& \Z[\U]_1\oplus \Z[\V]_1. 
\end{array}
\end{equation}

Our notation does not distinguish between $\delta_0$  (or $\delta_1$) in (\ref{sceq}) and 
the exponential mapping $\delta_0:K_0(\ac/\ec)\to K_1(\ec/\kc_\R)$ (or the index mapping
$\delta_1:K_1(\ac/\ec)\to K_0(\ec/\kc_\R)$, respectively) in the cyclic exact sequence of K-groups 
associated to (\ref{seq}). They are the same, modulo the isomorphism (\ref{k1ek}) (or (\ref{k0ek})). Still,
in the proofs of the following two propositions, we will have to be aware of how those isomorphisms are
defined. 

\begin{pro}In $($\ref{sceq}$)$, we have $\delta_0([[b(D)]_\ec]_0)=[\U]_1+[\V]_1=-\delta_0([[c(D)]_\ec]_0)$.
\label{del0}
\end{pro}

\pf Although $[c(D)]_{\kc_{\R}}$ is not a projection in $\ac/\kc_\R$, 
it is a self-adjoint lift for $[c(D)]_\ec$, since $c(D)$ is self-adjoint already in $\ac$. We therefore have 
\[
\delta_0([c(D)]_\ec)\,=\,-[e^{2\pi i[c(D)]_{\kc_{\R}}}]_1\,=\,-[[e^{2\pi ic(D)}]_{\kc_{\R}}]_1.
\]
In view of the isomorphism (\ref{e/k}), in order to write the above element of $K_1(\ec/\kc_\R)$ as a linear combination of 
$[\U]_1$ and $[\V]_1$ we must look at $\gamma_{_{T}}$, where $T=e^{2\pi i c(D)}=t(D)$, $t(\xi)=e^{2\pi ic(\xi)}$. By Theorem~\ref{th2},
\begin{equation}
\label{hoc}\gamma_{t(D)}(e^{2\pi i\varphi},+1)=
\gamma_{t(D)}(e^{2\pi i\varphi},-1)=Y_{\varphi}M_\varphi^tY_{-\varphi},\ \mbox{for all}\ \varphi\in\R.
\end{equation}

We have chosen $c$ so that 
\begin{equation}
t(\xi)=1  \ \ \mbox{if}\ \ \ |\xi|\geq 1/5.
\end{equation} 
This implies that
\begin{equation}
\label{viva}
t(j-\varphi)=1\ \ \mbox{if}\ \ 0\neq j\in\Z\ \ \mbox{and}\ \ |\varphi|\leq\frac{1}{2}.
\end{equation}
We now choose $E$ as the orthogonal projection onto the subspace of $\ell^2(\Z)$
generated by ${\mathbf d}=(d_n)_{n\in\Z}$, $d_n=0$ if $n\neq 0$, and $d_0=1$.
It follows from (\ref{viva}) that 
\begin{equation}
M_\varphi^t\ =\ I+[t(-\varphi)-1]E\ \ \mbox{if}\ \ |\varphi|\leq\frac{1}{2}.
\end{equation}

Let us choose a smooth $h:\R\to\R$ which vaninshes on $[-1/4,+1/4]$. For each $x\in[0,1]$ 
and each $\varphi\in[-1/2,+1/2]$, define
\[
U_x(\varphi)=Y_{(1-x)\varphi+xh(\varphi)}[I+(t(-\varphi)-1)E]Y_{(x-1)\varphi-xh(\varphi)}.
\]
$x\mapsto U_x$ is homotopy of unitaries in $C([-1/2,+1/2],\kc_\Z)^+$ satisfying 
\[
U_0(\varphi)=Y_{\varphi}M_\varphi^tY_{-\varphi}\ \ \mbox{and}\ \ 
U_1(\varphi)=Y_{h(\varphi)}[I+(t(-\varphi)-1)E]Y_{-h(\varphi)},
\]
for all $\varphi\in[-1/2,+1/2]$. 
Since $Y_{h(\varphi)}=I$ if $|\varphi|\leq 1/4$ and $t(-\varphi)-1=0$ if $|\varphi|\geq 1/5$,
we have $U_1(\varphi)=I+(t(-\varphi)-1)E$ for all $\varphi\in[-1/2,+1/2]$. Now 
$U_1(\varphi)=I$ for $|\varphi|\geq 1/5$ implies also that 
\[
V_1(z)=I+[g(z)-1]E,\ \  z\in\s^1,
\]
with $g(e^{2\pi i\varphi})=t(-\varphi)$ if $\varphi\in[-1/2,+1/2]$, defines a unitary 
$V_1\in C(\s^1,\kc_\Z)^+$ homotopic to $z\mapsto\gamma_{t(D)}(z,+1)=\gamma_{t(D)}(z,-1)$ (because of 
(\ref{hoc})). Since the winding number of $g$ around the origin is $+1$, $V_1$ is homotopic to 
$I+(\mbox{{\tt z}}-1)E$. This shows that 
\[
\delta_0([[c(D)]_\ec]_0)=-[(I+(\mbox{{\tt z}}-1)E,I+(\mbox{{\tt z}}-1)E)]_1=-([\U]_1+[\V]_1).
\]
Analogously, one can show that 
\[
\delta_0([[b(D)]_\ec]_0)=-[(I+(\overline{\mbox{{\tt z}}}-1)E,I+(\overline{\mbox{{\tt z}}}-1)E)]_1=[\U]_1+[\V]_1,
\]
as we wanted.
\cqd

Up to here we have understood $K_0$-elements as formal differences of classes of self-adjoint idempotents and 
$K_1$-elements as classes of unitaries. At this point, however, it is more convenient to use Banach-algebra K-theory 
definitions: $K_0$ will now consist of formal differences of classes of idempotents, and $K_1$ of classes of 
invertible elements. These two descriptions of the K-groups are equivalent for any \cst-algebra \cite{B}. Allowing $K_0$-elements 
to be represented by non-selfadjoint idempotents yields the very explicit expression for the index mapping given in the 
following lemma, which has been used before in \cite{N}, and is used below  in our computation (Proposition~\ref{ind}) of the 
index mapping $\delta_1$ in (\ref{sceq}).

\begin{lem} Let $A$ be a unital Banach algebra and let $\delta_1:K_1(A/J)\to K_0(J)$ denote the index mapping in the standard 
K-theory cyclic exact sequence associated to the short exact sequence of Banach algebras 
$0\to J\to A {\mathop{\to}\limits^{\pi}} A/J\to 0$.
If $u\in M_n(A/J)$ is an invertible, $\pi(a)=u$ and $\pi(b)=u^{-1}$, then 
\[
\delta_1([u]_1)=
\left[\left(\begin{array}{cc}
2ab-(ab)^2&a(2-ba)(1-ba)\\
(1-ba)b&(1-ba)^2
\end{array}\right)\right]_0-
\left[\left(\begin{array}{cc}
1&0\\0&0
\end{array}\right)\right]_0.
\]
\label{m}
\end{lem}
\pf Let $w\in M_{2n}(A)$ be defined by 
\[
w=\left(\begin{array}{cc}
2a-aba&ab-1\\1-ba&b
\end{array}\right).
\]
Then $w$ is invertible, 
\[
w^{-1}=\left(\begin{array}{cc}
b&1-ba\\ab-1&2a-aba
\end{array}\right)\ \ \mbox{and}\ \  \pi(w)=\left(\begin{array}{cc}u&0\\0&u^{-1}\end{array}\right).
\]
Our claim now follows immediately from \cite[8.3.1]{B}.\cqd

\begin{pro} Let $A\in\ac$ be such that $[A]_\ec$ is invertible in $\ac/\ec$. Then $\gamma_{_{A}}(z,-1)$ and  
$\gamma_{_{A}}(z,+1)$ are Fredholm operators in $\bc$ for every $z\in\s^1$ and 
\[
\delta_1([[A]_\ec]_1)\ =(\ind(\gamma_{_{A}}(1,-1))[E]_0,\ind(\gamma_{_{A}}(1,+1))[E]_0)\in K_0(\kc_\Z)\oplus K_0(\kc_\Z)
\]
\label{ind}
where $\ind$ denotes the Fredholm index, and $E$ is as in $($\ref{ktc}$\,)$.
\label{del1}
\end{pro}

\pf Let $B\in\ac$ be such that $I-AB$ and $I-BA$ belong to $\ec$. 
By Theorem~\ref{th3}, $\gamma_{_{I-AB}}=\gamma_{_{I}}-\gamma_{_{A}}\gamma_{_{B}}$
and $\gamma_{_{I-BA}}=\gamma_{_{I}}-\gamma_{_{B}}\gamma_{_{A}}$ belong to $C(S^*\s^1,\kc_\Z)$; i.e., for each $z\in\s^1$,
$\gamma_{_{B}}(z,\pm 1)$ is an inverse modulo $\kc_\Z$ for $\gamma_{_{A}}(z,\pm 1)$, what proves the first assertion.

Let us take $a=[A]_{\kc_{\R}}$ as a lift (pre-image) for $u=[A]_\ec$ with respect to the canonical mapping 
$\ac/\kc_\R\to\ac/\ec$. Let $B\in\ac$ be such that $[B]_\ec=u^{-1}$. Then $b=[B]_{\kc_{\R}}$ is a lift for $u^{-1}$.
Applying Lemma~\ref{m} and using the isomorphism (\ref{k0ek}), we see that $\delta_1([\U]_1)$ is equal to 
\begin{equation}
\label{mr}
\left(\ \left[\left(\begin{array}{cc}
2\mbox{{\tt A}}^{-}\mbox{{\tt B}}^{-}-(\mbox{{\tt A}}^{-}\mbox{{\tt B}}^{-})^2
&\mbox{{\tt A}}^{-}(2I-\mbox{{\tt B}}^{-}\mbox{{\tt A}}^{-})(I-\mbox{{\tt B}}^{-}\mbox{{\tt A}}^{-})\\
(I-\mbox{{\tt B}}^{-}\mbox{{\tt A}}^{-})\mbox{{\tt B}}^{-}&(I-\mbox{{\tt B}}^{-}\mbox{{\tt A}}^{-})^2
\end{array}\right)\right]_0-
\left[\left(\begin{array}{cc}
I&0\\0&0
\end{array}\right)\right]_0\right.,
\end{equation}
\[
\left.\left[\left(\begin{array}{cc}
2\mbox{{\tt A}}^{+}\mbox{{\tt B}}^{+}-(\mbox{{\tt A}}^{+}\mbox{{\tt B}}^{+})^2
&\mbox{{\tt A}}^{+}(2I-\mbox{{\tt B}}^{+}\mbox{{\tt A}}^{+})(I-\mbox{{\tt B}}^{+}\mbox{{\tt A}}^{+})\\
(I-\mbox{{\tt B}}^{+}\mbox{{\tt A}}^{+})\mbox{{\tt B}}^{+}&(I-\mbox{{\tt B}}^{+}\mbox{{\tt A}}^{+})^2
\end{array}\right)\right]_0-
\left[\left(\begin{array}{cc}
I&0\\0&0
\end{array}\right)\right]_0\ \right),
\]
where $\mbox{{\tt A}}^{\pm}=\gamma_{_{A}}(1,\pm 1)$ and $\mbox{{\tt B}}^{\pm}=\gamma_{_{B}}(1,\pm 1)$.

Let us recall that the index mapping $\delta_F$ for the exact sequence
\[
0\longrightarrow\kc_\Z\longrightarrow\bc\longrightarrow\bc/\kc_\Z\longrightarrow 0
\]
``is'' the Fredholm index \cite[9.4.2]{Ro}, in the sense that $\delta_F$  maps $[[T]_{\kc_{\Z}}]_1$, if $T\in\bc$ is a 
Fredholm operator, into $\ind(T)[E]_0$. On the other hand, again by Lemma~\ref{m},
the first component in the expression for $\delta_1([\U]_1)$ in (\ref{mr}) equals 
$\delta_F([[\mbox{{\tt A}}^{-}]_{\kc_{\Z}}]_1)$. 
I.e., that first component is equal to $\ind(\mbox{{\tt A}}^{-})[E]_0$, as we wanted.
Analogously for the other component. \cqd

\begin{cor} In $($\ref{sceq}$)$, we have $\delta_1([[A_1]_\ec]_1)=(0,[E]_0)$, 
$\delta_1([[A_2]_\ec]_1)=([E]_0,0)$, $\delta_1([[A_3]_\ec]_1)=(0,-[E]_0)$ and $\delta_1([[A_4]_\ec]_1)=(-[E]_0,0)$. 
\label{cor2}
\end{cor}

\pf It takes a straighforward application of all definitions to write down explicit formulas for  
$\gamma_{_{A_{i}}}(1,\pm 1)$, $i=1,\cdots, 4$. Their indices 
can be computed by \cite[19.1.14]{H}, for example. We omit the details, which can be found in \cite{S}.\cqd

We are now ready to conclude the analysis of (\ref{sceq}). It follows 
from Corollary~\ref{cor2} that $\delta_1$ is surjective, and hence the upper right arrow 
in (\ref{sceq}), which maps $[[I]_{\kc_{\R}}]_0$ to $[[I]_{\ec}]_0$, is injective. 
By Corollary~\ref{cor1}, $[b(D)]_{\ec}$ and $[c(D)]_{\ec}$ are mutually orthogonal projections in $\ac/\ec$, hence:
\[
[[b(D)]_{\ec}]_0+[[c(D)]_{\ec}]_0= [[b(D)+c(D)]_{\ec}]_0=[[I]_{\ec}]_0.
\]
By Proposition~\ref{del0}, the kernel of $\delta_0$ is generated by
$[[b(D)]_{\ec}]_0+[[c(D)]_{\ec}]_0$, and this shows that 
\[
K_0(\ac/\kc_\R)=\Z[[I]_{\kc_{\R}}]_0.
\]

The image of $\delta_0$ is complemented by $\Z[\U]_1$ (by Proposition~\ref{del0}) and the kernel of $\delta_1$ is equal 
to $\Z([[A_1]_\ec]_1+[[A_3]_\ec]_1)\oplus\Z([[A_2]_\ec]_1+[[A_4]_\ec]_1)$ (by Corollary~\ref{cor2}). The exactness of (\ref{sceq})
then implies that 
\[
K_1(\ac/\kc_\R)=\Z[[B_1]_{\kc_{\R}}]_1\oplus\Z[[B_2]_{\kc_{\R}}]_1\oplus\Z[[B_3]_{\kc_{\R}}]_1,
\]
if $B_1\in\ac$, $B_2\in\ac$ and $B_3\in\ec^+$ are such that 
\begin{equation}
\label{b1b2}
\sigma_{_{B_{1}}}=\sigma_{_{A_{1}}}\cdot\sigma_{_{A_{3}}}=(l,l), \ \ 
\sigma_{_{B_{2}}}=\sigma_{_{A_{2}}}\cdot\sigma_{_{A_{4}}}=(\tilde l,\tilde l)
\end{equation}
and 
\begin{equation}
\label{b3}
[\gamma_{_{B_{3}}}]_1=[\U]_1. 
\end{equation}
If we take $B_1=l(M)$ and $B_2=\tilde l(M)$, then we obviously get (\ref{b1b2}). 
We claim that
\[
B_3=b(M)t(D)+c(M),
\]
$t$ as in the proof of Proposition~\ref{del0}, satisfies (\ref{b3}). That follows from 
\[
\gamma_{_{B_{3}}}(e^{2\pi i\varphi},+1)=Y_{\varphi}M_\varphi^tY_{-\varphi},\ \ 
\gamma_{_{B_{3}}}(e^{2\pi i\varphi},-1)=I
\]
and the fact, shown in the proof of Proposition~\ref{del0}, that the unitary
\[
e^{2\pi i\varphi}\mapsto Y_{\varphi}M_\varphi^tY_{-\varphi}
\] 
is homotopic, within $C(\s^1,\kc_\Z)^+$, to $I+(\mbox{{\tt z}}-1)E$. This settles the question of computing the 
K-theory of $\ac/\kc_\R$. 

The Fredholm index of $B_3$ is $(-1)$, by Theorem~\ref{thm4}. 
The index mapping in the six-term exact sequence associated to $0\to\kc_\R\to\ac\to\ac/\kc_\R\to 0$ is therefore surjective. 
This implies that the canonical projection $\ac\to\ac/\kc_\R$ induces isomorphisms between $K_0(\ac)$ and $K_0(\ac/\kc_\R)$
and between $K_1(\ac)$ and the kernel of the index mapping $K_1(\ac/\kc_\R)\to\Z$, which is the subgroup generated by 
$[[B_1]_{\kc_{\R}}]_1$ and $[[B_2]_{\kc_{\R}}]_1$ (since $B_1$ and $B_2$ are invertible). This finishes the proof of:

\begin{thm}\label{main}
$K_0(\ac)=\Z[I]_0$ and $K_1(\ac)=\Z[l(M)]_1\oplus \Z[\tilde l(M)]_1$.
\end{thm}

%
%
%
%
%
%

\section{The operator-valued symbol exact sequence}\label{img}

Let $\ac^\diamond$ denote the \cst-algebra of bounded operators on $L^2(\R)$ generated by all
operators of type (2) or (3) (listed at the beginning of this paper). In this section, we 
define a surjective \cst-algebra homomorphism 
\[
\psi:\ac\longrightarrow\ac^\diamond\oplus\ac^\diamond
\]
and exhibit a \cst-algebra isomorphism between the images of $\gamma$ and $\psi$
\[
\iota:\text{Im}\,\gamma\longrightarrow \text{Im}\,\psi
\]
such that $\psi=\iota\circ\gamma$. We compute the K-theory of $\ac^\diamond$, and then analyse the standard 
cyclic exact sequence associated to $\psi$ and $\gamma$. 

Given $u\in L^{2}(\R)$, for almost every $\varphi\in\R$, the sequence
$
u^\diamond(\varphi)=(u(\varphi-j))_{j\in\Z}
$ 
belongs to $\ell^2(\Z)$. A unitary mapping
\[
W:L^{2}(\R)\longrightarrow L^2(\s^1;\ell^2(\Z))
\]
is defined by
$(Wu)(e^{2\pi i\varphi})=Y_\varphi u^\diamond(\varphi)$. $Wu$ is smooth whenever $u$ is.
Given a continuous function $g\in C(\s^1,\bc)$, let $\mu(g)$ denote the bounded operator 
on $L^2(\s^1;\ell^2(\Z))$ of multiplication by $g$. An injective 
\cst-algebra homomorphism is defined by
\[
\begin{array}{rccc}
\mu\oplus \mu: &C(S^*\s^1,\bc)&\longrightarrow&{\mathcal L}(L^2(\s^1;\ell^2(\Z)))\oplus{\mathcal L}(L^2(\s^1;\ell^2(\Z)))\\
&(f,g)&\longmapsto&(\mu(f),\mu(g))
\end{array}
\]
(we have used notation explained just before Proposition \ref{prop1}).
It is straightforward to check that, given $A\in\ac^\diamond$ and $d\in C[-\infty,+\infty]$, we have 
\[
\mu\oplus\mu(\gamma_{_{B}})=(d(-\infty)WF^{-1}AFW^{-1},d(+\infty)WF^{-1}AFW^{-1}),
\] 
if $B=d(M)A$. 
The set of all such $B$'s and the compact ideal together generate a dense subalgebra of $\ac$. 
Since the compacts are contained in the kernel of $\gamma$, we get: 

\begin{thm}\label{au}
There exists a surjective \cst-algebra homomorphism
\[
\psi:\ac\longrightarrow \ac^\diamond\oplus\ac^\diamond
\]
such that $\psi(A)=(A,A)$ if $A\in\ac^\diamond$, $\psi(a(M))=(a(-\infty)I,a(+\infty)I)$ if
$a\in C([-\infty,+\infty])$ and $\psi(K)=0$ if $K$ is compact. Denoting by $\iota$ the conjugation by 
$FW^{-1}\oplus FW^{-1}$ composed with $\mu\oplus\mu$, we have $\iota\circ\gamma=\psi$.
\end{thm}

Thus $\img$ is isomorphic to $\ac^\diamond\oplus\ac^\diamond$.
We now compute the K-theory of $\ac^\diamond$.

\begin{thm}\label{ktad}
$K_0(\ac^\diamond)=\Z[I]_0$ and $K_1(\ac^\diamond)=\Z[e^{iM}]_1$, where $e^{iM}$ denotes the operator of
multiplication by $x\mapsto e^{ix}$. 
\end{thm}

\pf Let $\alpha$ denote the automorphism of $A=C[-\infty,+\infty]$ defined by 
\[
[\alpha(f)](x)=f(x-1),\ x\in \R,\ f\in A,
\] 
and let $A\times_{\alpha}\Z$ denote the envelopping \cst-algebra \cite{D} of the Banach algebra
with involution $\ell^1(\Z,A)$ of all summable $\Z$-sequences in $A$ equipped with the product
\[
({\mathbf f}\cdot{\mathbf g})(n)\ =\ \sum_{k\in\Z}f_k\alpha^k(g_{n-k}),\ \ 
n\in\Z,\ {\mathbf f}=(f_k)_{k\in\Z},\ {\mathbf g}=(g_k)_{k\in\Z},
\]
and involution ${\mathbf f}^*(n)=\alpha^{n}(\bar f_{-n})$. 

The K-groups of $A$ and $A\times_{\alpha}\!\Z$ fit \cite{PV} into the 
{\em Pimsner-Voiculescu exact sequence}
\begin{equation}
\label{pv}
\def\mapup#1{\Big\uparrow\rlap{$\vcenter{\hbox{$\scriptstyle#1$}}$}}
\def\mapdn#1{\Big\downarrow\rlap{$\vcenter{\hbox{$\scriptstyle#1$}}$}}
\begin{array}{ccccc}
K_0(A)&{\mathop{\longrightarrow}\limits^{\text{\tt id}-\alpha_{*}^{-1}}}&K_0(A)
&{\mathop{\longrightarrow}\limits^{i_*}}&K_0(A\times_{\alpha}\!\Z)
\\&&&&
\\\mapup{\delta_1}&&&&\mapdn{\delta_0}
\\&&&&
\\K_1(A\times_{\alpha}\!\Z)&{\mathop{\longleftarrow}\limits^{i_*}}&K_1(A)&
{\mathop{\longleftarrow}\limits^{\text{\tt id}-\alpha_{*}^{-1}}}&K_1(A),
\end{array}
\end{equation}
where $\text{\tt id}$ denotes the identity on K-groups and $i$ the canonical inclusion
of $A$ into $A\times_{\alpha}\!\Z$: $i(a)$ is the sequence whose only nonvanishing 
entry is $a$ at position zero. Abusing notation, we will denote $i(a)$ by $a$.
Using that $K_0(A)=\Z[1]_0$, that $K_1(A)=0$, and that $\alpha(1)=1$,
it follows easily from (\ref{pv}) that $K_0(A\times_{\alpha}\!\Z)=\Z[1]_0$ and that 
$K_1(A\times_{\alpha}\!\Z)$ is also isomorphic to $\Z$. 

For each integer $k$, let ${\mathbf d}_k$ denote the sequence whose only nonvanishing entry is $1$ 
at position $k$. ${\mathbf d}_1$ is a unitary in $A\times_{\alpha}\!\Z$ such that, for every $a\in A$, 
$\alpha(a)={\mathbf d}_1a{\mathbf d}_1^*$. It then follows from the equation at 
the beginning of \cite[page 102]{PV} that $K_1(A\times_{\alpha}\!\Z)=\Z[{\mathbf d}_{-1}]_1$.

To prove the theorem, it therefore suffices to show that there exists a \cst-algebra
isomorphism $\varphi:A\times_{\alpha}\!\Z\to\ac^\diamond$ such that $\varphi({\mathbf d}_{-1})=e^{iM}$.

For sequences ${\mathbf a}=(a_j)_{j\in\Z}$ with only a finite number of nonzero entries, define 
\[
\varphi({\mathbf a})=\sum_{j}a_j(D)e^{-ijM}.
\]
It is elementary to check that $\varphi$ extends to a *-homomorphism from $\ell^1(\Z,A)$ to 
$\ac^\diamond$, which obviously satisfies $\varphi({\mathbf d}_{-1})=e^{iM}$. It then follows from \cite[2.7.4]{D} 
that $\varphi$ extends to a \cst-algebra homomorphism $\varphi:A\times_{\alpha}\!\Z\to\ac^\diamond$. The image of $\varphi$ 
is dense, by definition of $\ac^\diamond$. So, all that remains to be shown is that $\varphi$ is injective.

Given $z\in\s^1$ and ${\mathbf a}=(a_j)_{j\in\Z}\in\ell^1(\Z,A)$, define $\gamma_z({\mathbf a})=(z^ja_j)_{j\in\Z}$.
$\gamma_ z$ extends to an automorphism of $A\times_{\alpha}\!\Z$ and $z\mapsto\gamma_z$ is an action of the circle. 
For each $z\in\s^1$, let $U_z$ denote the unitary mapping on $L^2(\R)$ given by $(U_zf)(x)=z^xf(x)$, $x\in\R$, and define
$\beta_z(A)=FU_zF^{-1}AFU_z^{-1}F^{-1}$, $A\in\ac^\diamond$. $z\mapsto\beta_z$ is an action of $\s^1$ by automorphisms on
$\ac^\diamond$ satisfying $\beta_z\circ\varphi=\varphi\circ\gamma_z$ for each $z$. To show that $\varphi$ is injective, it 
therefore suffices to show that its restriction to the fixed points of $\gamma$ is injective (by \cite[Proposition 2.9]{E}).

Let $C$ denote the algebra of points fixed by $\gamma$,
\[
C=\{{\mathbf x}\in A\times_{\alpha}\!\Z;\,\gamma_ z({\mathbf x})={\mathbf x}\ \text{for all}\ z\in\s^1\},
\]
and let $E:A\times_{\alpha}\!\Z\to A\times_{\alpha}\!\Z$ be defined by
\[
E({\mathbf x})=\int_{\s^1}\gamma_z({\mathbf x})\slash\!\!\!d\!z,\ \slash\!\!\!d(e^{i\theta})=d\theta/2\pi.
\]
If ${\mathbf x}\in C$, then $E({\mathbf x})={\mathbf x}$; hence $C$ is contained in the image of $E$, $\text{\tt Im}\,E$. 
Since, for all $j\neq 0$ the integral $\int z^j\slash\!\!\!d\!z$ vanishes, we get $E({\mathbf a})\in A$ for 
all ${\mathbf a}\in\ell^1(\Z,A)$. Since $A$ is a closed subalgebra of $A\times_{\alpha}\!\Z$ and $\ell^1(\Z,A)$
is dense in $A\times_{\alpha}\!\Z$, we get $\text{Im}\,E\subseteq A$. Since $A\subseteq C$, we get $A=C$. 
But it is clear that $\varphi$ restricted to $A$ is injective, as we wanted. \cqd

Let $\jc$ denote the \cst-algebra of bounded operators on $L^2(\R)$ generated by all operators 
of the form $a(M)b(D)$, for $a\in C_0(\R)$ and $b\in C[-\infty,+\infty]$. Choose a sequence $\chi_j\in C_0(\R)$
such that $0\leq\chi_j\leq 1$, the support of each $\chi_j$ is contained in the interval $(-j,+j)$ and 
$\chi_j(x)\to 1$ for all $x\in\R$. Using the facts stated in the Introduction about the commutators $[a(M),b(D)]$, it is 
straightforward to prove (this is a particular case of \cite[Lemma
VII.1.2]{C2}) that a given $A\in\ac$ belongs to $\jc$ if and only if there exists a sequence of 
compact operators $C_j$ such that $\chi_ j(M)A+C_j\to A$. 
This implies that $\gamma$ vanishes on a dense subset of $\jc$; and hence $\jc\subset\ker\gamma$. 
It also shows that $A\in\jc$ if and only if $\sigma_{_{A}}$ vanishes over all points of the form $((\pm\infty,z),\pm\infty)$,
$z\in\s^1$. Since the supremum of $|\sigma_{_{A}}|$ over all those points is bounded by $||\gamma_{_{A}}||$ 
\cite[Proposition 3.4]{CM}, we then get 
\begin{equation}
\label{kerg}
\jc=\ker\gamma.
\end{equation} 

The restriction of $\sigma$ to $\jc$ induces an isomorphism
\begin{equation}
\label{mca}
\jc/\kc_\R\cong C_0(\R\times\{-\infty,+\infty\}).
\end{equation} 
This follows from a general result \cite[Theorem VI.2.2]{C2} about comparison
algebras ($\jc$ is the {\em minimal comparison algebra} over $\R$).
A direct proof for our case is given in \cite[Lema 1.15]{S}. 

If we now quotient the exact sequence $0\to\ker\gamma\to\ac\to\img\to 0$ by the compacts and use Theorem~\ref{au} and (\ref{kerg}),
we obtain the exact sequence of \cst-algebras
\begin{equation}
\label{esg}
0\longrightarrow\jc/\kc_\R\,{\mathop{\longrightarrow}\limits^{\iota}}\,\ac/\kc_\R\,
{\mathop{\longrightarrow}\limits^{\psi}}\,\ac^\diamond\oplus\ac^\diamond\longrightarrow 0.
\end{equation}
and its associated six-term exact sequence 
\begin{equation}
\label{csg}
\def\mapdn#1{\Big\downarrow\rlap{$\vcenter{\hbox{$\scriptstyle#1$}}$}}
\begin{array}{ccccc}
0&\longrightarrow&K_0(\ac/\kc_\R)&{\mathop{\longrightarrow}\limits^{\psi_*}}&K_0(\ac^\diamond\oplus\ac^\diamond)
\\
&&&&
\\
\Big\uparrow&&&&\mapdn{\delta_0}
\\
&&&&
\\
K_1(\ac^\diamond\oplus\ac^\diamond)&{\mathop{\longleftarrow}\limits^{\psi_*}}
&K_1(\ac/\kc_\R)&{\mathop{\longleftarrow}\limits^{\iota_*}}& K_1(\jc/\kc_\R).
\end{array}
\end{equation}

It follows from Theorem~\ref{ktad} that
$K_0(\ac^\diamond\oplus\ac^\diamond)=\Z[(I,0)]_0\oplus\Z[(0,I)]_0$. Since
$\psi(I)=(I,I)$, we have
\begin{equation}
\label{kerdz}
[(I,0)]_0+[(0,I)]_0\ \in\ \ker\delta_0.
\end{equation}
The isomorphism (\ref{mca}) implies that $[[U_+]_{\kc_\R}]_1$ and $[[U_-]_{\kc_\R}]_1$ generate
$K_1(\jc/\kc_\R)\cong\Z\oplus\Z$, if each $U_\pm$ is a Fredholm operator in $\jc\oplus\C I$ such that its
principal symbol is equal to one on $\R\times\{\mp\infty\}$ and winds once around the origin over 
$\R\times\{\pm\infty\}$. We may take, for example, $U_+=e^{2\pi ib(M)}b(D)+c(D)$ and
$U_-=e^{2\pi ic(M)}c(D)+b(D)$ ($b$ and $c$ were defined at the beginning of Section~\ref{principal}).
Since the winding numbers of $\bar\sigma(U_+)$ and $\bar\sigma(U_+)$ are equal, it follows from 
Theorem~\ref{thm4} that 
\begin{equation}
\label{ccs}
i_*([[U_+]_{\kc_\R}]_1)=i_*([[U_-]_{\kc_\R}]_1),
\end{equation} 
where $i:\jc/\kc_\R\to\Cc/\kc_\R$ denotes the canonical inclusion. From (\ref{ccs}), we then get 
\begin{equation}
\label{iota}
\iota_*([[U_+]_{\kc_\R}]_1)=\iota_*([[U_-]_{\kc_\R}]_1)
\end{equation}
($\iota:\jc/\kc_\R\to\ac/\kc_\R$ denotes the inclusion, as defined in (\ref{esg})).

The exactness of (\ref{csg}), together with (\ref{kerdz}) and (\ref{iota}), yields  
\[
K_0(\ac/\kc_\R)\,=\,\Z[[I]_{\kc_\R}]_0
\]
and the short exact sequence
\begin{equation}
\label{ses}
0\longrightarrow\Z[[U_+]_{\kc_\R}]_1 {\mathop{\ \longrightarrow\ }\limits^{\iota_*}}  K_1(\ac/\kc_\R)
{\mathop{\ \longrightarrow\ }\limits^{\psi_*}}
K_1(\ac^\diamond\oplus\ac^\diamond)\longrightarrow 0.
\end{equation}
Since $\psi(l(M))=(I,e^{iM})$ and $\psi(\tilde l(M))=(e^{iM},I)$, Theorem~\ref{ktad} and (\ref{ses}) 
together imply 
\[ 
K_1(\ac/\kc_\R)\ =\ \Z[[U_+]_{\kc_\R}]_1\oplus\Z[[l(M)]_{\kc_\R}]_1\oplus\Z[[\tilde l(M)]_{\kc_\R}]_1.
\]
Since, by Theorem \ref{thm4}, the Fredholm index of $U_+$ is (-1), this gives another proof of Theorem~\ref{main}.
\section*{Acknowledgements}
We are grateful to Ruy Exel for many helpful conversations.
It was his the idea of using that $\varphi$ is a 
covariant homomorphism in the proof of Theorem~\ref{ktad}.

Severino Melo was partially supported by CNPq (Processo 306214/2003-2), and C\'{\i}ntia Silva had
a scholarship from CAPES during her graduate studies.

{\footnotesize 

\vskip2.0cm
 
Instituto de  Matem\'atica e Estat\'{\i}stica, 
Universidade de S\~ao Paulo

Caixa Postal 66281,
05311-970 S\~ao Paulo, Brazil.

toscano@ime.usp.br, cintia@ime.usp.br}


\begin{thebibliography}{99}

\bibitem{AS} {\sc M. F. Atiyah \&\ I. M. Singer}, {\em The index of elliptic operators I}, Annals of Math. 
{\bf 87} (1968), 484-530.

\bibitem{B} {\sc B. Blackadar}, K-theory for operator algebras, Cambridge
University Press, Cambridge, 1998. 

\bibitem{C1} {\sc H. O. Cordes}, Elliptic pseudodifferential operators -
an abstract theory, Lecture Notes in Mathematics {\bf 756}, Springer, 
Berlin, 1979.

\bibitem{C2} {\sc H. O. Cordes}, Spectral theory of linear differential 
operators and comparison algebras, London Mathematical Society Lecture 
Note Series {\bf 76}, Cambridge University Press, Cambridge, 1987. 

\bibitem{CP} {\sc H. O. Cordes}, {\em On the two-fold symbol chain of a C$^*$-algebra of 
singular integral operators on a polycilinder}, Revista Mat. Iberoamericana {\bf 2} (1986), 215-234.

\bibitem{CH} {\sc H. O. Cordes \&\ E. A. Herman}, {\em Gelfand theory of 
pseudodifferential operators}, Amer. J. Math {\bf 90} (1968), 681-717.

\bibitem{CM} {\sc H. O. Cordes \&\ S. T. Melo}, {\em An algebra of singular
integral operators with kernels of bounded oscillation and application to 
periodic differential operators}, J. Differential Equations {\bf 75}-2 (1988),
216-238.

\bibitem{D} {\sc J. Dixmier}, Les C$^*$-alg\`ebres et leurs repr\'esentations,
Gauthier-Villars, Paris, 1969.

\bibitem{E} {\sc R. Exel}, {\em Circle actions on C$^*$-algebras, partial 
automorphisms and a generalized Pimsner-Voiculescu exact sequence}, 
J. Funct. Anal. {\bf 122}-2 (1994), 361-401.

\bibitem{H} {\sc L.H\"ormander}, The analysis of linear partial differential 
operators III, Springer, New York, 1990.

\bibitem{Lt} {\sc R. Lauter}, Holomorphic functional calculus in several variables 
and $\Psi^*$-algebras of totally characteristic operators on manifolds with boundary, 
Thesis, Mainz, 1996.

\bibitem{L} {\sc R. Lauter}, {\em Pseudodifferential analysis on conformally 
compact spaces}, Mem. Amer. Math. Soc. {\bf 163}-777, 2003.

\bibitem{M1} {\sc S. T. Melo}, Comparison Algebras with Periodic Symbols,
Thesis, Berkeley, 1988. 

\bibitem{M2} {\sc S. T. Melo}, {\em A comparison algebra on a cylinder with 
semi-periodic multiplications}, Pacific J. Math. {\bf 146}-2 (1990), 281-304.

\bibitem{MNS} {\sc S. T. Melo, R. Nest \&\ E. Schrohe}. {\em $C\sp *$-structure and 
$K$-theory of Boutet de Monvel's algebra}, J. Reine Angew. Math. {\bf 561} (2003),
145-175.

\bibitem{MSS}{\sc S. T. Melo, T. Schick \&\ E. Schrohe, }{\em A K-theoretic proof of 
Boutet de Monvel's index theorem}, Arxiv, math.KT/0403059, to appear in  
J. Reine Angew. Math.

\bibitem{Me} {\sc R. B. Melrose}, {\em The eta invariant and families of pseudodifferential 
operators}, Math. Res. Lett. {\bf 2} (1995), 541-561.

\bibitem{MoN} {\sc R. B. Melrose \&\ V. Nistor}, {\em $K$-theory of $C\sp *$-algebras 
of $b$-pseudodifferential operators}, Geom. Funct. Anal. {\bf 8}-1 (1998), 88-122.

\bibitem{MN} {\sc B. Monthubert \&\ V. Nistor}, {\em A topological index theorem for 
manifolds with corners}, Arxiv, math.KT/0507601.

\bibitem{Mo} {\sc S. Moroianu}, {\em K-Theory of suspended pseudo-differential operators},
K-Theory {\bf 28} (2003), 167-181.

\bibitem{N} {\sc V. Nistor}, {\em Higher index theorems and the boundary map in
cyclic cohomology}, Doc. Math. {\bf 2} (1997), electronic, 263-295. 

\bibitem{PV}  {\sc M. Pimsner \&\ D. Voiculescu}, {\em Exact sequences for K-groups
and Ext-groups of certain cross-product C$^*$-algebras}, 
J. Operator Theory {\bf 4}-1 (1980), 93-118. 

\bibitem{P} {\sc S. C. Power}, {\em Commutator ideals and pseudodifferential 
\cst-algebras}, Quart. J. Math. Oxford  Ser. (2) {\bf 31}-124 (1980), 467-489.

\bibitem{R} {\sc V. S. Rabinovich}, {\em On the algebra generated by the 
pseudodifferential operators on $R\sp{n}$, the operators of multiplication 
by almost periodic functions and the shift operators},  
Soviet Math. Dokl. {\bf 25}-2 (1982), 498-502.

\bibitem{Ro} {\sc M. R\o rdam, F. Larsen \&\ N. Laustsen}, An Introduction to K-theory
for \cst-algebras, Cambridge University Press, Cambridge, 2000.

\bibitem{S} {\sc C. C. Silva}, K-teoria de operadores pseudodiferenciais
com s\'{\i}mbolos semiperi\'odicos, Tese, S\~ao Paulo (Arxiv, math.OA/0505261), 2005.

\end{thebibliography}
\end{document}